\documentclass{article}

\usepackage{arxiv}

\usepackage[utf8]{inputenc} 
\usepackage[T1]{fontenc}    
\usepackage{hyperref}       
\usepackage{url}            
\usepackage{booktabs}       
\usepackage{amsfonts}       
\usepackage{nicefrac}       
\usepackage{microtype}      
\usepackage{lipsum}		
\usepackage{amsmath}
\usepackage{enumerate}
\usepackage{tikz}
\usepackage{algorithm}
\usepackage{algorithmic}
\usepackage{booktabs}
\newtheorem{definition}{Definition}
\newtheorem{theorem}{Theorem}
\newtheorem{lemma}{Lemma}

\newtheorem{proof}{Proof}
\def\QEDopen{{\setlength{\fboxsep}{0pt}\setlength{\fboxrule}{0.2pt}\fbox{\rule[0pt]{0pt}{1.3ex}\rule[0pt]{1.3ex}{0pt}}}}
\def\QED{\QEDopen}
\def\proof{\noindent{\bf Proof}: }
\def\endproof{\hspace*{\fill}~\QED\par\endtrivlist\unskip}
\newtheorem{remark}{Remark}
\DeclareSymbolFont{ugmL}{OMX}{mdugm}{m}{n}
\SetSymbolFont{ugmL}{bold}{OMX}{mdugm}{b}{n}
\DeclareMathAccent{\wideparen}{\mathord}{ugmL}{"F3}

\title{A Combinatorial Algorithm for the Multi-commodity Flow Problem}

\author{
  Pengfei Liu\\
  Building No.1, Zhonguancun Road\\
  Haidian District\\
  Beijing, China \\
  \texttt{pengfeil89@foxmail.com} \\
}

\begin{document}
\maketitle

\begin{abstract}
This paper researches combinatorial algorithm for the multi-commodity flow problem.   We relax the capacity constraints and introduce a \emph{penalty function} \(h\) for each arc. If the flow exceeds the capacity on arc \(a\), arc \(a\) would have a penalty cost. Based on the \emph{penalty function} \(h\), a new conception , \emph{equilibrium pseudo-flow}, is introduced. If  the equilibrium pseudo-flow is  a nonzero-equilibrium pseudo-flow, there exists no feasible solution for the multi-commodity flow problem; if  the equilibrium pseudo-flow is  a zero-equilibrium pseudo-flow, there exists  feasible solution for the multi-commodity flow problem and the zero-equilibrium pseudo-flow is the feasible solution. Then   a \emph{non-linear} description of the multi-commodity flow problem is given, whose solution is equilibrium pseudo-flow. And  a combinatorial algorithm(Frank-Wolfe algorithm) is designed to obtain equilibrium pseudo-flow.  Besides, the content in this paper can be easily generalized to minimum cost multi-commodity flow problem.
\end{abstract}


\section{Introduction}
The \emph{multi-commodity flow problem} (MFP) is the problem of designing flow of several different commodities through a common network with arc capacities.
Given a directed graph \(G(V,A)\), a capacity function \(u : A \rightarrow Q^+\), \(K\) origin-destination pairs of nodes , defined by \(K_k=(s_k, t_k, d_k)\)  where \(s_k\) and \(t_k\) are
 the origin and destination of commodity \(k\), and \(d_k\) is the demand. The flow of commodity \(k\) along arc \((i,j)\) is \(f^k_{ij}\). The objective is to obtain an assignment of flow which
 satisfies the demand for each commodity without violating the capacity constraints. The constraints can be summarized as follows:

 \begin{equation}\label{eq1}
\begin{split}
 &\sum_{k\in K}f^k_{ij}\leq u_{ij}, \forall (i,j)\in A\\
  &\sum_{j\in \delta^+(i)}f^k_{ij}-\sum_{j\in \delta^-(i)}f^k_{ji}=\left\{
  \begin{aligned}
  &d_k, \text{ if } i=s_k\\
  &-d_k, \text{ if } i=t_k\\
  &0,\text{ if } i\in V - \{s_k,t_k\}\\
  \end{aligned}
  \right.\\
  & f^k_{ij}\geq 0, \forall k\in K, (i,j) \in A\\
\end{split}
\end{equation}
where $\delta^+(i)=\{j|(i,j)\in A\}, \delta^-(i)=\{j|(j,i)\in A\}$. In this paper, we assume $s_k\neq t_k$. The first expression is capacity constraint. The second is flow conservation constraint and the last is non-negative constraint.

Multicommodity flow problems have attracted great attention since the publication of
the works of \cite{ff} and \cite{hu}. \cite{assad} gives a comprehensive survey, which includes decomposition, partitioning, compact inverse methods, and primal-dual algorithms. Although there are many combinatorial algorithms for single-commodity flow models like Ford-Fulkerson algorithm(\cite{ff}), Edmonds-Karp algorithm(\cite{Edmonds}), Dinic's algorithm (\cite{Dinic}) and push-relabel algorithm(\cite{Goldberg}), there is no known combinatorial algorithm for multi-commodity flow problem.
It  is well known that  MFP can be solved in polynomial time using linear programming. However, up to date,
there is no other way to solve the problem precisely without using linear programming. In this paper, we would \emph{give the first combinatorial algorithm} for the multi-commodity flow problem.

\subsection{Our contribution}
A new conception, \emph{equilibrium pseudo-flow}, is introduced and a \emph{non-linear} description of the multi-commodity flow problem is given whose solution is equilibrium pseudo-flow. Besides, we give a combinatorial algorithm to obtain  equilibrium pseudo-flow for the multi-commodity flow problem. To the best of our knowledge, this is the first algorithm to obtain the solution of  the multi-commodity flow problem without using linear programming.

\section{Equilibrium Pseudo-flow}
\label{sec:headings}

Unlike other methods, our combinatorial algorithm \emph{does not maintain the capacity constraints} throughout the execution. The algorithm, however, maintains a \emph{pseudo-flow}, which
is a function \(\mathbf{f} : K \times V \times V \rightarrow \mathbb{R}^+\) that just satisfies the flow conservation on every node. That is, a pseudo-flow \(\mathbf{f}\) is a feasible solution of Expression~(\ref{eq2}).
\begin{equation}\label{eq2}
\begin{split}
&\sum_{j\in \delta^+(i)}f^k_{ij}-\sum_{j\in \delta^-(i)}f^k_{ji}=\left\{
  \begin{aligned}
  &d_k, \text{ if } i=s_k\\
  &-d_k, \text{ if } i=t_k\\
  &0,\text{ if } i\in V - \{s_k,t_k\}\\
  \end{aligned}
  \right.\\
  & f^k_{ij}\geq 0, \forall k\in K, (i,j) \in A\\
\end{split}
\end{equation}

We introduce a \emph{penalty function} \(h\) for each arc \((i,j)\in A\), which is defined as
\begin{equation}\label{eq3}
h(f_{ij})=\left\{
\begin{aligned}
0 ~~~~~~~~~~~~~~& if~f_{ij} \leq u_{ij} \\
f_{ij} - u_{ij}~~~~~ & if~f_{ij} > u_{ij}
\end{aligned}
\right.
\end{equation}

If the flow on an arc \((i,j)\) is less than the capacity, the penalty of arc \((i,j)\) is zero. Otherwise, the penalty of arc \((i,j)\)  is the amount by which the
flow  exceeds the capacity. 

Intuitively, the greater the \(h(f_{ij})\) is, the more 'congested' the arc \((i,j)\) is.  By using \(\{h(f_{ij}), \forall (i,j) \in A\}\)  as weights for the arcs, the longer the path \(p_{s_kt_k}\)
is, the  more 'congested' the path \(p_{s_kt_k}\) is, where \(p_{s_kt_k}\) is a path connecting  \(s_k\) and \(t_k\). In fact, for a pair\((s_k,t_k)\), our algorithm iteratively adjusts the flow to the shortest paths until all the used paths have equal length.

We introduce the concept of \emph{equilibrium pseudo-flow} here, which is the key to the combinatorial
algorithm.

\begin{definition}
\label{def1}
By using \{ \(h(f_{ij}), \forall (i,j) \in A\)\}  as weights for all the arcs, a pseudo-flow \(\mathbf{f}\) is called an \emph{equilibrium pseudo-flow} if it satisfies the following conditions:
\begin{enumerate}[(i)]
\item for any given pair \((s_k,t_k)\), all used paths connecting \(s_k\) and \(t_k \) have equal and minimum length;
\item for any given pair \((s_k,t_k)\), all unused paths connecting \(s_k\) and \(t_k \) have greater or equal length;
\end{enumerate}
\end{definition}
where a path \(p\) connecting \(s_k\) and \(t_k \) is called \emph{used} if there exists \(s_k-t_k\) flow on path \(p\), otherwise it is called \emph{unused}. The conditions above are also called \emph{equilibrium conditions}. Note that the conception above is similar to 'user equilibrium'(\cite{Wardrop}),  which is a sound and simple behavioral principle to describe the spreading of trips.

\begin{definition}
\label{def2}
 An equilibrium pseudo-flow \(\mathbf{f}\) is called  zero-equilibrium pseudo-flow if \{ \(h(f_{ij})=0, \forall (i,j) \in A\)\}. Otherwise, it is called nonzero-equilibrium pseudo-flow.
\end{definition}

Obviously, by the definition above, a zero-equilibrium pseudo-flow is a feasible flow that satisfies Expression~(\ref{eq1}). Therefore, we have the following theorem:
\begin{theorem}
\label{theo1}
Given \(\{(s_k,t_k,d_k): k \in K\}\) and capacity reservation \(\{u_{ij}: (i,j)\in A\}\), the feasible region of Expression~(\ref{eq1}) is not empty if and only if there exists a zero-equilibrium pseudo-flow.
\end{theorem}

In fact, if there exists a nonzero-equilibrium pseudo-flow, there is no feasible solution for Expression~(\ref{eq1}). Before proving this conclusion, we need the following lemma, which  was originally given  by \cite{f} and  \cite{b},  and subsequently observed by \cite{d}.

\begin{lemma}
\label{lem1}
Given \(\{(s_k,t_k,d_k): k \in K\}\) and capacity reservation \(\{u_{ij}: (i,j)\in A\}\), the feasible region of Expression~(\ref{eq1}) is not empty if and only if:
\begin{equation}\label{eq11}
\sum_{k \in K}l_{s_k,t_k}^\mu d_k\leq \sum_{(i,j)\in A}\mu_{ij}u_{ij}, \forall\mu: A\rightarrow \mathbb{Z}^+\cup \{0\}\\
\end{equation}
where \(l_{s_k,t_k}^\mu\) is the length of the shortest path from \(s_k\) to \(t_k\) using \(\mu\) as weights for the arcs.
\end{lemma}
\begin{theorem}
\label{theo2}
Given \(\{(s_k,t_k,d_k): k \in K\}\) and capacity reservation \(\{u_{ij}: (i,j)\in A\}\), the feasible region of Expression~(\ref{eq1}) is empty if there exists a nonzero-equilibrium pseudo-flow.
\end{theorem}
\proof
Let \(f^p_k\) be the flow on path \(p\) connecting \(s_k\) and \(t_k \)  and \(\delta_{a,p}^k\) indicator variable where
\begin{equation}\label{eq4}
\delta_{a,p}^k=\left\{
\begin{aligned}
1~~~~& \text{if arc \(a\) is on path \(p\) connecting  \(s_k\) and \(t_k \) } \\
0~~~~~&  \text{otherwise}
\end{aligned}
\right.
\end{equation}
Let \(P_k\) be the set of all the used paths connecting \(s_k\) and \(t_k \) , we have
\begin{equation}\label{eq5}
\sum_{p\in P_k}f^p_k = d_k~~~\forall k\in K
\end{equation}
Let  \(l_{s_k,t_k}\) be the length of the shortest path from \(s_k\) to \(t_k\) and \(l_{s_k,t_k}^p\) the length of the  path \(p\) connecting \(s_k\) and \(t_k\) using the penalty function \(\{h(f_a): \forall a \in A\}\) as weights for the arcs. The following formulation shows the relationship between \(l_{s_k,t_k}^p\) and \(\{h(f_a): \forall a \in A\}\).

\begin{equation}\label{eq6}
l_{s_k,t_k}^p = \sum_{a\in A} h(f_a) \delta_{a,p}^k
\end{equation}
Based on the relationship between arc flows and path flows, the following equation holds:

\begin{equation}\label{eq7}
f_a = \sum_{k\in K}  \sum_{p \in P_k} \delta_{a,p}^k f_p^k
\end{equation}
According to the definition of the equilibrium pseudo-flow,  all used paths connecting \(s_k\) and \(t_k \)  have equal and minimum length, that is,
\begin{equation}\label{eq30}
\left\{
\begin{aligned}
l_{s_k,t_k}^p = l_{s_k,t_k} \text{~if~} f_k^p > 0\\
l_{s_k,t_k}^p \geq l_{s_k,t_k} \text{~if~} f_k^p = 0\\
\end{aligned}
\right.
\end{equation}

 Then we have
\begin{equation}\label{eq8}
\begin{split}
\sum_{k\in K} l_{s_k,t_k}d_k &= \sum_{k\in K} l_{s_k,t_k}(\sum_{p \in P_k}f_p^k)~~~~~~~~~~~~~~~~~~\setminus\setminus \text{  by Expression~(\ref{eq5})}\\
&= \sum_{k\in K}  (\sum_{p \in P_k}l_{s_k,t_k} f_p^k)\\
&= \sum_{k\in K}  (\sum_{p \in P_k}l_{s_k,t_k}^pf_p^k)~~~~~~~~~~~~~~~~~~~\setminus\setminus \text{  by Expression~(\ref{eq30})}\\
&= \sum_{k\in K}  (\sum_{p \in P_k}(\sum_{a\in A} h(f_a) \delta_{a,p}^k) f_p^k)~~~~~~~~\setminus\setminus \text{  by Expression~(\ref{eq6})}\\
&= \sum_{a\in A} h(f_a)(\sum_{k \in K}  \sum_{p \in P_k} \delta_{a,p}^k f_p^k)\\
&= \sum_{a\in A} h(f_a)f_a~~~~~~~~~~~~~~~~~~~~~~~~~~~\setminus\setminus \text{  by Expression~(\ref{eq7})}
\end{split}
\end{equation}

 According to the definition of the nonzero-equilibrium pseudo-flow, there exists at least an arc \(a\) that satisfies \(f_a > u_a\). Since \(h(f_a)\) = 0 if \(f_a\leq u_a\) and  \(h(f_a)> 0 \) if \(f_a > u_a\),
  \(\sum_{a\in A} h(f_a)f_a > \sum_{a\in A} h(f_a)u_a\). That is,

\begin{equation}\label{eq9}
\begin{split}
\sum_{k \in K} l_{s_k,t_k} d_k &= \sum_{a\in A} h(f_a)f_a~~~~~~~~~~~~~~~~\setminus\setminus \text{  by Formulation~(\ref{eq8})}\\
& > \sum_{a\in A} h(f_a)u_a
\end{split}
\end{equation}

By Lemma~\ref{lem1}(viewing \(h\) as \(\mu\)),the feasible region of Expression~(\ref{eq1}) is empty.
\endproof

\begin{remark}
In fact, Theorem~\ref{theo2} is a necessary  and sufficient condition. We would see that in  Section~\ref{s4}.
\end{remark}

Assume we have an algorithm to get the equilibrium pseudo-flow. Based on Theorem~\ref{theo1} and Theorem~\ref{theo2}, we have the following conclusion:

\begin{theorem}
\label{theo3}
If  the equilibrium pseudo-flow is  a nonzero-equilibrium pseudo-flow, there exists no feasible solution for Expression~(\ref{eq1}); if  the equilibrium pseudo-flow is  a zero-equilibrium pseudo-flow, there exists  feasible solution for Expression~(\ref{eq1}) and the zero-equilibrium pseudo-flow is a feasible solution.
\end{theorem}

So what we need to do is only to design an algorithm to obtain the equilibrium pseudo-flow.

\section{The Formulation of MFP}
\label{s4}
 The multi-commodity flow problem (MFP) is always regarded as a linear programming problem. However, in this part, we will give a non-linear programming formulation of  MFP, whose solution is an equilibrium pseudo-flow.

\subsection{The Basic Formulation}

Let \(f_a\) be the sum of the flow of all pairs on arc \(a\) and \(h(f_a)\) be penalty   function on arc \(a\).
\begin{equation}\label{eq12}
\begin{split}
{\bf min}\quad &z = \sum_a\int_0^{f_a}h(\omega)d\omega\\
{\bf s.t}\quad &\sum_{j\in \delta^+(i)}f^k_{ij}-\sum_{j\in \delta^-(i)}f^k_{ji}=\left\{
  \begin{aligned}
  &d_k, \text{ if } i=s_k\\
  &-d_k, \text{ if } i=t_k\\
  &0,\text{ if } i\in V - \{s_k,t_k\}\\
  \end{aligned}
  \right.\\
  & f^k_{ij}\geq 0, \forall k\in K, (i,j) \in A\\
\end{split}
\end{equation}

In the program above, the objective function is the sum of the integrals of
the arc penalty function. The first  constraint is flow conservation constraint and the second is non-negative constraint. Note that there is no capacity constraint here. According to the definition of the penalty function \(h\), if the feasible region of Expression (1) is not empty, \emph{the minimum value of the objective function is zero; otherwise it is greater than zero.}

The formulation above is similar to Beckmann Formulation (\cite{Beckmann}), whose solution is called User Equilibrium	(\cite{Wardrop}). However, \cite{Beckmann} didn't give an reasonable interpretation of the objective function. It is just viewed strictly as a mathematical construct that is utilized to solve User Equilibrium problems. In this paper we give an economic interpretation of the objective function.

Let's look at a simple example. Assume there are ten cars queueing up to cross an intersection. The intersection allows a car to pass at one time and each  car will take 1 unit time to go through the intersection. Obviously, after 10 units time all the cars would go through the intersection. Now let's look this phenomenon from another perspective. The time the \(ith\) car spends to go through the intersection is \(i\) units time because it needs to wait until the cars  in front go through the intersection. Therefore, the sum of the time of every car to go through the intersection is \(1+2+3+\cdots+10 = 55\).  That is, the sum of the time of every car to go through the intersection is 55 and the time of the last car to go through the intersection is 10. Now let's look at the objective function.  The penalty of that the \(ith\) unit flow passes through the arc \(a\) is \(h(i)\). The integration \(\int_0^{f_a}h(\omega)d\omega\) means the sum of the penalty of every unit flow to pass through the arc \(a\). \emph{For an arc \(a\), what the objective function minimizes is the sum of the penalty of every unit flow to pass through the arc \(a\), not the penalty of the last unit flow to pass through the arc \(a\).}

\subsection{Equivalence}

To demonstrate the equivalence between the equilibrium pseudo-flow and
Program ~(\ref{eq12}), it has to be shown that any flow pattern that solves Program ~(\ref{eq12})
satisfies the equilibrium conditions. This equivalency is demonstrated in this
part by proving that the Karush-Kuhn-Tucker conditions for Program ~(\ref{eq12}) are identical to the equilibrium conditions.

\begin{lemma}
\label{lemma2}
\(t(f_a) = \int_0^{f_a}h(\omega)d\omega\) is a convex function.
\end{lemma}
\proof{Proof.}
The derivative  of \(t(f_a)\) is \(h(f_a)\), which is monotone nondecreasing function. So \(t(f_a) = \int_0^{f_a}h(\omega)d\omega\) is a convex function.
\endproof

\begin{lemma}
\label{lemma3}
Let \(\mathbf{f}^\ast\) be a solution of Program ~(\ref{eq12}).  \(\mathbf{f}^\ast\) is the optimal solution of Program ~(\ref{eq12}) if and only if  \(\mathbf{f}^\ast\) satisfies the  Karush-Kuhn-Tucker conditions of  Program ~(\ref{eq12}).
\end{lemma}
\proof{Proof.}
By Lemma~\ref{lemma2}, \(t(f_a) = \int_0^{f_a}h(\omega)d\omega\) is a convex function. Therefore, the objective function \(z = \sum_a\int_0^{f_a}h(\omega)d\omega\) is a convex function. Besides, the inequality constraints of  Program ~(\ref{eq12}) are continuously differentiable convcave functions and the equality constraints of Program ~(\ref{eq12}) are affine functions. So Karush-Kuhn-Tucker conditions are necessary and sufficient for optimality of Program~(\ref{eq12}) (\cite{Boyd}).
\endproof

Obviously, Program~(\ref{eq12}) is a minimization problem with nonnegativity constraints and linear equality. The  Karush-Kuhn-Tucker conditions of such formulation are as following:

\begin{equation}\label{eq13}
\begin{split}
&{\bf Stationarity}\\
& \quad -\frac{\partial z}{\partial f^k_{ij}} = -\mu_{ij}^k + (\lambda_i^k - \lambda_j^k),\forall k\in K, (i,j) \in A\\
&{\bf Primal~feasibility}\\
& \quad \sum_{j\in \delta^+(i)}f^k_{ij}-\sum_{j\in \delta^-(i)}f^k_{ji}=\left\{
  \begin{aligned}
  &d_k, \text{ if } i=s_k\\
  &-d_k, \text{ if } i=t_k\\
  &0,\text{ if } i\in V - \{s_k,t_k\}\\
  \end{aligned}
  \right.\\
& \quad -f^k_{ij}\leq 0, \forall k\in K, (i,j) \in A\\
&{\bf Dual~feasibility}\\
& \quad \mu_{ij}^k \geq 0,\forall k\in K, (i,j) \in A\\
&{\bf Complementary~slackness}\\
& \quad \mu_{ij}^kf_{ij}^k = 0,\forall k\in K, (i,j) \in A\\
\end{split}
\end{equation}
Obviously,
\begin{equation*}
\frac{\partial z}{\partial f_{ij}^k}  = h(f_{ij}) \frac{\partial f_{ij}}{\partial f_{ij}^k} = h(f_{ij})
\end{equation*}
Substituting the  expression above into Stationarity expression in KKT conditions,
\begin{equation*}
  h(f_{ij}) = \mu_{ij}^k + ( \lambda_j^k - \lambda_i^k ),\forall k\in K, (i,j) \in A\\
\end{equation*}
For a path \(p = (s_k, v_1, v_2, \cdots, v_m, t_k)\), the length of \(p\) is
\begin{equation*}
\begin{split}
  length(p) &= h(f_{s_kv_1}) + h(f_{v_1v_2}) + \cdots + h(f_{v_{m-1}v_m}) + h(f_{v_mt_k})\\
  &= \mu_{s_kv_1}^k + (\lambda_{v_1}^k - \lambda_{s_k}^k) +  \mu_{v_1v_2}^k + (\lambda_{v_2}^k - \lambda_{v_1}^k) + \cdots + \mu_{v_mt_k}^k +(\lambda_{t_k}^k - \lambda_{v_m}^k)\\
  &= \lambda_{t_k}^k - \lambda_{s_k}^k + \mu_{s_kv_1}^k + \mu_{v_1v_2}^k + \cdots + \mu_{v_mt_k}^k
\end{split}
\end{equation*}

The condition above holds for every path between any  pair in the network. For an \(arc~(i,j)\) on a used path \(p_{used}\) between \(pair~(s_k,t_k)\), the flow \(f_{ij}^k\) is greater than zero. By complementary slackness \( \mu_{ij}^kf_{ij}^k = 0\) in KKT conditions, we have \( \mu_{ij}^k = 0\). Therefore,
\begin{equation*}
\begin{split}
  length(p_{used}) &= \lambda_{t_k}^k - \lambda_{s_k}^k + \mu_{s_kv_1}^k + \mu_{v_1v_2}^k + \cdots + \mu_{v_mt_k}^k\\
   &= \lambda_{t_k}^k - \lambda_{s_k}^k
\end{split}
\end{equation*}
By the expression above,  all the used paths between \(pair~(s_k,t_k)\) have the same length \((\lambda_{t_k}^k - \lambda_{s_k}^k)\).

For an unused path \(p_{unused}\) between \(pair~(s_k,t_k)\), the length of \(p_{unused}\) is
\begin{equation*}
  length(p_{unused}) = \lambda_{t_k}^k - \lambda_{s_k}^k + \mu_{s_kv_1}^k + \mu_{v_1v_2}^k + \cdots + \mu_{v_mt_k}^k
\end{equation*}
By dual feasibility \( \mu_{ij}^k \geq 0 \) in KKT conditions, \(length(p_{unused})\) is greater or equal to \(length(p_{used})\).

With this interpretation above, it is now clear that:
\begin{enumerate}[(i)]
\item all the used paths connecting   \(s_k\)and \(t_k \) have equal and minimum length;
\item all the unused paths connecting \(s_k\)and \(t_k \) have greater or equal length;
\end{enumerate}
That is, the optimal solution of Program~(\ref{eq12}) is an equilibrium pseudo-flow.

\subsection{Frank-Wolfe Algorithm}

The Program~(\ref{eq12})  includes a convex  objective function, a linear constraint set and a non-negative constraint set, which could be efficiently solved by Frank-Wolfe algorithm (\cite{Wolfe}). Applying Frank-Wolfe algorithm to Program~(\ref{eq12}), at the \(nth\) iteration, needs the following linear program:

\begin{equation}\label{eq15}
\begin{split}
{\bf min}\quad &z^n(\mathbf{y}) = \sum_{k,ij}\frac{\partial z(\mathbf{f}_n)}{\partial f_{ij}^k}y^k_{ij} = \sum_{k,ij}h(f_{ij,n})y^k_{ij}\\
{\bf s.t}\quad &\sum_{j\in \delta^+(i)}y^k_{ij}-\sum_{j\in \delta^-(i)}y^k_{ji}=\left\{
  \begin{aligned}
  &d_k, \text{ if } i=s_k\\
  &-d_k, \text{ if } i=t_k\\
  &0,\text{ if } i\in V - \{s_k,t_k\}\\
  \end{aligned}
  \right.\\
  & y^k_{ij}\geq 0, \forall k\in K, (i,j) \in A\\
\end{split}
\end{equation}
where \(f_{ij,n}\) is the flow on arc \((i,j)\)  at the \(nth\) iteration.

Note that this program doesn't have capacity constraints and the penalties  are not flow-dependent. In other words, the  program minimizes the total penalties over a network with fixed
penalties \(\{h(f_{ij,n}): \forall (i,j) \in A\}\). Obviously, the penalties  will be minimized by assigning all \(s_k-t_k\) flows to the shortest path connecting \(s_k\) and \(t_k\).  Such an assignment is performed by computing the shortest paths between all pairs.  Since the penalty of each arc is 0 at \(0th\) iteration, we can establish the initial pseudo-flow \(\mathbf{f}\) by assigning all the \(s_k-t_k\) flow to certain path connecting \(s_k\) and \(t_k\). Therefore, The Frank-Wolfe algorithm applied to solve Program~(\ref{eq12}) can be given as follows:
\begin{algorithm}
\caption{ Frank-Wolfe Algorithm applied to MCF}
\label{alg:B}
\begin{algorithmic}
\STATE ~~~
\STATE \emph{Initialization}: examine each pair \((s_k, t_k)\) in turn and assign all the \(s_k-t_k\) flow to certain path connecting \(s_k\) and \(t_k\). This yields \(\mathbf{f}_1\). Set \(n = 1\);
\REPEAT
\STATE  \emph{Update}: set \( \{h(f_{ij,n}): \forall (i,j) \in A\}\) as the weights of every arc;
\STATE \emph{ Direction-finding}: compute the shortest paths between all pairs and assigning all \(s_k-t_k\) flows to the shortest path connecting \(s_k\) and \(t_k\), which yields \(\mathbf{y}_n\).
\STATE \emph{ Line search}: find \(\alpha_n\) by solving \({\bf min}_{0\leq\alpha \leq 1} \sum_{(i,j)\in A}\int_0^{f_{ij,n} + \alpha(y_{ij,n} - f_{ij,n}) }h(\omega)d\omega\)
\STATE \emph{Move}: set \(f_{ij,n+1} = f_{ij,n} + \alpha_n(y_{ij,n} - f_{ij,n})\)
\UNTIL{some convergence criterion is met }
\end{algorithmic}
\end{algorithm}

\begin{remark}
The convergence of the Frank–Wolfe algorithm above is sublinear, that is, the error in the objective function to the optimum is \(O(1/k)\) after \(k\) iterations\cite{Wolfe}.
\end{remark}

\begin{remark}
In fact, all the conclusion in this paper is true if the penalty function \(h\)  satisfies the following definition,
\begin{equation*}
h(f_{ij})=\left\{
\begin{aligned}
0 ~~~~~~~~~~~~~~& if~f_{ij} \leq u_{ij} \\
g(f_{ij} - u_{ij})~~~~~ & if~f_{ij} > u_{ij}
\end{aligned}
\right.
\end{equation*}
where \(g(0) = 0\) and \(g(x)\) is strictly monotone increasing function when \(x \geq 0\).
\end{remark}

\begin{remark}
If the penalty function \(h\) is defined as following,
\begin{equation*}
h(f_{ij})=\left\{
\begin{aligned}
c_{ij} ~~~~~~~~~~~~~~& if~f_{ij} \leq u_{ij} \\
c_{ij} + M(f_{ij} - u_{ij})~~~~~ & if~f_{ij} > u_{ij}
\end{aligned}
\right.
\end{equation*}
Program~(\ref{eq12}) is a description of \emph{minimum cost multi-commodity flow problem}, where \(M\) is big enough and \(\{c_{ij}:(i,j)\in A\}\) is the cost of every arc. Therefore, by defining the penalty function \(h\) as above, the content in this paper can be easily generalized to minimum cost multi-commodity flow problem.
\end{remark}

\begin{remark}
There exist  several variants of the Frank-Wolfe algorithm which have faster convergence(\cite{Jaggi,Pena}). We would not discuss here.
\end{remark}

\begin{remark}
Since the objective function \(z\) is convex, we have \(z(\mathbf{x}) \geq z(\mathbf{y}) + (\mathbf{x - y})^T\nabla z(\mathbf{y})\), which  also holds for the  optimal solution \(\mathbf{f^*}\). That is, \(z(\mathbf{f^*}) \geq z(\mathbf{f^n}) + (\mathbf{f^*} - \mathbf{f^n})^T\nabla z(\mathbf{f^n}) \geq min_{\mathbf{x} \in D} \{z(\mathbf{f^n}) + (\mathbf{x - f^n})^T\nabla z(\mathbf{f^n})\} = z(\mathbf{f^n}) - \mathbf{f^n}^T\nabla z(\mathbf{f^n}) + min_{\mathbf{x} \in D} \{\mathbf{x} ^T\nabla z(\mathbf{f^n})\}\), where \(\mathbf{f^n}\) is the flow  at the \(nth\) iteration and \(D\) is the feasible region.
If \( z(\mathbf{f^n}) - \mathbf{f^n}^T\nabla z(\mathbf{f^n}) + min_{\mathbf{x} \in D} \{\mathbf{x} ^T\nabla z(\mathbf{f^n})\} > 0\), the objective function \(z\) must be greater than zero and there exists no zero-equilibrium pseudo-flow. The algorithm can be terminated in advance.
\end{remark}

\section{Conclusion}
This paper gives a combinatorial algorithm for the multi-commodity flow problem.  Unlike other methods, the combinatorial algorithm does not maintain the capacity constraints throughout the execution. The algorithm, however, maintains a pseudo-flow, which just satisfies the flow conservation on every nodes.  We introduce a \emph{penalty function} \(h\) for each arc, which is positively related to the quantity that the flow exceeds the capacity. Then  a \emph{non-linear} description of the multi-commodity flow problem is given whose solution is equilibrium pseudo-flow and we also give a combinatorial algorithm  to obtain equilibrium pseudo-flow. Generally, we give a new idea to solve the multi-commodity flow problem, but it need further studies  to design more effective algorithm  to obtain the equilibrium pseudo-flow.



\end{document}